\documentclass[11pt]{article}
\usepackage{mathrsfs}
\usepackage{CJK}
\usepackage{amsmath}
\usepackage{fancyhdr}

\begin{document}
\title{ A Note On The Jacobian Conjecture}
\author{{ Dan Yan\footnote{E-mail address:yan-dan-hi@163.com}}    \\
\small School of Mathematical Sciences, Graduate University of
Chinese Academy of Sciences,\\
\small Beijing 100049, China\\}
\date{}
\maketitle \begin{abstract} In this note, we show that, if the Druzkowski mappings $F(X)=X+(AX)^{*3}$, i.e. $F(X)=(x_1+(a_{11}x_1+\cdots+a_{1n}x_n)^3,\cdots,x_n+(a_{n1}x_1+\cdots+a_{nn}x_n)^3)$, satisfies $TrJ((AX)^{*3})=0$, then $rank(A)\leq \frac{1}{2}(n+\delta)$ where $\delta$ is the number of diagonal elements of A which are equal to zero. Furthermore, we show the Jacobian Conjecture is true for the Druzkowski mappings in dimension $\leq 9$ in the case $\prod_{i=1}^{n}a_{ii}\neq0$.
\end{abstract} {\bf Keywords.} Jacobian conjecture,
Polynomial mapping, Druzkowski mapping\\
{\bf MSC(2010).} Primary 14E05  Secondary 14A05;14R15 \vskip 2.5mm

Let
$F=(F_1(x_1,\cdots,x_n),\cdots,F_n(x_1,\cdots,x_n))^t:\bf{C}^n\rightarrow\bf{C}^n$
be a polynomial mapping, that is,
$F_i(x_1,\cdots,x_n)\in\bf{C}[x_1,\cdots,x_n]$ for all $1\leq i\leq
n$. Let $JF=(\frac{\partial F_i}{\partial x_j})_{n\times n}$ be the
Jacobian matrix of $F$. The well-known Jacobian Conjecture(JC)
raised by O.H. Keller in 1939 (\cite{1}) states that a polynomial
mapping $F:\bf{C}^n\rightarrow\bf{C}^n$ is invertible if the
Jacobian determinant $|JF|$ is a nonzero constant. This conjecture
has being attacked by many people from various research fields and
remains open even when $n=2$! (Of course, a positive answer is
obvious when $n=1$ by the elements of linear algebra.)  See \cite{2}
and \cite{3} and the references therein for a wonderful 70-years
history of this famous conjecture. It can be easily seen that JC
is true if JC holds for all polynomial mappings whose Jacobian
determinant is 1. We make use of this convention in the present
paper. Among the vast interesting and valid results, a relatively
satisfactory result obtained by S.S.S.Wang(\cite{4}) in 1980 is that
JC holds for all polynomial mappings of degree 2 in all dimensions.
The most powerful and surprising result is the reduction to degree
3, due to H.Bass, E.Connell and D.Wright( \cite{2}) in 1982 and
A.Yagzhev(\cite{5}) in 1980, which asserts that JC is true if JC
holds for all polynomial mappings of degree 3(what is more, if JC
holds for all cubic homogeneous polynomial mappings!). In the same spirit of the above degree reduction method, another efficient way to tackle JC is the Druzkowski's Theorem(\cite{6}): JC is true if it is true for all Druzkowski mappings (in all dimension $\geq 2$).

\indent Recall that $F$ is a cubic homogeneous map if $F=X+H$ with
$X$ the identity (written as a column vector) and each component of
$H$ being either zero or cubic homogeneous. A cubic homogeneous
mapping $F=X+H$ is a {\it\bf Druzkowski (or cubic linear) mapping}
if each component of $H$ is either zero or a third power of a linear
form. Each Druzkowski mapping $F$ is associated to a scalar matrix
$A$ such that $F=X+(AX)^{*3}$, where $(AX)^{*3}$ is the {\it\bf
Druzkowski symbol} for the vector $(A^1X)^3,\cdots, (A^nX)^3)$ with
$A^i$ the $i-$th row of $A$. Clearly, a Druzkowski mapping is
uniquely determined by this matrix $A$.\\
\indent

{\bf Theorem } {\it Let $F=X+H$ be a Druzkowski mappings. If $TrJ((AX)^{*3})=0$, then $rank(A)\leq \frac{1}{2}(n+\delta)$ where $\delta$ is the number diagonal elements of A which are equal to zero.}\\
\indent

{\bf Proof of Theorem } Set $t_i=a_{i1}x_1+a_{i2}x_2+\cdots+a_{in}x_n$ for $1\leq i\leq n$. Since $TrJ((AX)^{*3})=0$, therefore, we have $$a_{11}t_1^2+a_{22}t_2^2+\cdots+a_{nn}t_n^2=0$$
i.e. $$(x_1,x_2,\cdots,x_n)[a_{11}\left(\begin {array}{c}a_{11}\\a_{12}\\\vdots\\a_{1n}\end{array}\right)(a_{11},a_{12},\cdots,a_{1n})+\cdots+a_{nn}\left(\begin {array}{c}a_{n1}\\a_{n2}\\\vdots\\a_{nn}\end{array}\right)(a_{n1},a_{n2},\cdots,a_{nn})]\left(\begin {array}{c}x_{1}\\x_{2}\\\vdots\\x_{n}\end{array}\right)=0$$
Since $(A^tDA)^t=A^tDA$, therefore, we have $A^tDA=0$ where\\
\[D=\left[\begin{array}{cccc}
a_{11}&0&\cdots&0\\
0&a_{22}&\cdots&0\\
\cdots&\cdots&\cdots&\cdots\\
0&0&\cdots&a_{nn}\end{array}\right].\]\\
Set $r(A)=r$. Then there exists invertible matrices $P$ and $Q$ such that $PAQ=E_r$ and $(PAQ)^t=E_r$ where\\
\[E_r=\left[\begin{array}{cccccc}
1&0&\cdots&0&\cdots&0\\
0&1&\cdots&0&\cdots&0\\
\cdots&\cdots&\cdots&\cdots&\cdots&\cdots\\
0&0&\cdots&1&\cdots&0\\
0&0&\cdots&0&\cdots&0\\
\cdots&\cdots&\cdots&\cdots&\cdots&\cdots\\
0&0&\cdots&0&\cdots&0\end{array}\right].\]\\
Therefore, we have $A^tDA=(Q^t)^{-1}E_r(P^t)^{-1}DP^{-1}E_rQ^{-1}=0$. Set $H=(P^t)^{-1}DP^{-1}$ and $H=(h_{ij})_{n\times n}$. Since $P$ is invertible, so we have $rank(H)=rank(D)=n-\delta$. Therefore, we get $E_rHE_r=0$\\
i.e. \\
\[\left[\begin{array}{ccccccc}
h_{11}&h_{12}&\cdots&h_{1r}&0&\cdots&0\\
h_{21}&h_{22}&\cdots&h_{2r}&0&\cdots&0\\
\cdots&\cdots&\cdots&\cdots&\cdots&\cdots&\cdots\\
h_{r1}&h_{r2}&\cdots&h_{rr}&0&\cdots&0\\
0&0&\cdots&0&0&\cdots&0\\
\cdots&\cdots&\cdots&\cdots&\cdots&\cdots&\cdots\\
0&0&\cdots&0&0&\cdots&0\end{array}\right]=0.\]\\
Since $rank(H)=n-\delta$, so $r\leq n+\delta-r$. Therefore, we have $r\leq \frac{1}{2}(n+\delta)$, which completes the proof.\\

Since the Jacobian Conjecture is true for all cubic homogeneous polynomials in dimension $\leq 4$ (see [8]) and we know if $rankA\leq 4$, $F$ and $G'$ are are paired, where $G'$ are cubic homogeneous polynomials in dimension no more than 4 which related to F.(see Theorem2 in [9]). Therefore, we have the following conclusion:\\
\indent

{\bf Corollary } {\it Let $F=X+H$ be a Druzkowski mapping, if $\prod_{i=1}^{n}a_{ii}\neq0$, then the Jacobian Conjecture is true for the Druzkowski mappings in dimension $\leq 9$.}\\
\indent
\indent

{\bf Proof of Corollary } since $detJF=1$, we have $TrJ((AX)^{*3})=0$. if $\prod_{i=1}^{n}a_{ii}\neq0$, we have $\delta=0$. Therefore, $rank(A)\leq \frac{n}{2}$. If $n\leq9$, we have $rank(A)\leq4$. then F and $G'$ are Gorni-Zampieri pairings, where $G'$ are 4 dimension cubic homogeneous polynomials and $detJG'=1$. Since in this situation, the Jacobian Conjecture is true, so we have the conclusion.\\

\indent
\indent
In the proof of the theorem, we have $r(A)\leq \frac{n}{2}$. Now we give an example to show that is the best bound for $rank(A)$ in such case. \\
\indent
 {\bf Example } {\it Let $n=4$, and}
    \begin{equation}
  \left\{ \begin{aligned}
  F_1=(x_1+ix_2+x_3+x_4)^3+x_1\\
  F_2=i(x_1+ix_2+x_3+x_4)^3+x_2\\
  F_3=-(x_1+ix_2-x_3+x_4)^3+x_3\\
  F_4=-(x_1+ix_2-x_3+x_4)^3+x_4.\\
                          \end{aligned} \right.
  \end{equation}
 Then we have \\
\[A=\left[\begin{array}{cccc}
1&i&1&1\\
-i&1&-i&-i\\
-1&-i&1&-1\\
-1&-i&1&-1\\\end{array}\right],\]\\
and clearly, $F$ is invertible and $rank(A)=2$.

  Acknowledgment: The author is very grateful to Professor Yuehui Zhang who introduced the Conjecture and gave great help when the author studied the problem and Professor Guoping Tang who read the paper carefully and gave some good advice.

\end{document}